\documentclass[11pt,twoside]{amsart}

\usepackage{amssymb, palatino,graphicx, hyperref}
\usepackage[normalem]{ulem}
\usepackage[all]{xy}
\usepackage{xcolor}
\usepackage{ifthen}

\graphicspath{{./}{./Figures/}}

\newtheorem{dfn}{Definition}[section]
\newtheorem{thm}[dfn]{Theorem}
\newtheorem{lem}[dfn]{Lemma}
\newtheorem{rem}[dfn]{Remark}

\newenvironment{preuve}{{\em \bf Proof:}}{\hfill $\blacksquare$}

\def\P{\ensuremath{\mathbb P}}

\def\R{\ensuremath{\mathbb R}}
\def\Z{\ensuremath{\mathbb Z}}
\def\H{\ensuremath{\mathbb H}}
\def\C{\ensuremath{\mathbb C}}

\def\CP1{\ensuremath{\mathbb C \mathbb P^1}}

\def\Pn-1{\ensuremath{\P^{n-1}}}

\definecolor{red}{rgb}{.6,0,0}
\definecolor{green}{rgb}{0,.6,0}
\definecolor{darkgreen}{rgb}{0,0.3,0}
\definecolor{purple}{rgb}{0.5,0,0.5}
\definecolor{darkblue}{rgb}{0,0,0.7}
\definecolor{greenblue}{rgb}{0,0.4,0.5}
\definecolor{myblue}{HTML}{1685d6}
\definecolor{mypurple}{HTML}{b82e6c}

\newboolean{draft}

\newcommand{\cmt}[1]
{\ifthenelse {\boolean{draft}}
{{\sc \tiny \color{red} #1}}
{}}

\newcommand{\newbb}[1]
{\ifthenelse {\boolean{draft}}
{{\color{darkblue} #1}}
{#1}}

\newcommand{\newbbb}[1]
{\ifthenelse {\boolean{draft}}
{{\color{greenblue} #1}}
{#1}}

\newcommand{\nopost}[1]
{\ifthenelse {\boolean{draft}}
{{\color{cyan} #1}}
{}}

\newcommand{\maynopost}[1]
{\ifthenelse {\boolean{draft}}
{{\color{purple} #1}}
{}}

\newcommand{\margincmt}[1]
{\ifthenelse {\boolean{draft}}
{\marginpar{{\sc \tiny \color{red} #1}}}
{}}
\newcommand{\inred}[1]
{\ifthenelse{\boolean{draft}}{{\color{red} #1}}{#1}}

\newcommand{\new}[1]
{\ifthenelse {\boolean{draft}}
{{\color{green} #1}}
{#1}}

\newcommand{\neww}[1]
{\ifthenelse {\boolean{draft}}
{{\color{darkgreen} #1}}
{#1}}

\newcommand{\newb}[1]
{\ifthenelse {\boolean{draft}}
{{\color{blue} #1}}
{#1}}

\newcommand{\del}[1]
{\ifthenelse {\boolean{draft}}
{{\color{magenta} #1}}
{}}

\newboolean{details_on}

\newcommand{\details}[1]
{\ifthenelse {\boolean{details_on}}
{{\color{darkgreen} \tiny #1}}
{}}

\title{Generalized Laurent monomials in nonrational toric geometry}
\author{Fiammetta Battaglia}
\author{Elisa Prato}
\thanks{The authors were partially supported by the PRIN Project ``Real and Complex Manifolds: Topology, Geometry and Holomorphic Dynamics" (MIUR, Italy) and by GNSAGA (INdAM, Italy).}

\setboolean{draft}{true}
\setboolean{details_on}{true}

\begin{document}

\begin{abstract} We generalize Laurent monomials to toric quasifolds, a special class of
highly singular spaces that extend simplicial toric varieties to the nonrational setting. 
\end{abstract}
\maketitle
\begin{small}

\medskip
\noindent \textbf{Mathematics~Subject~Classification:}
14M25, 52B20, 53D20.
\end{small}

\section*{Introduction}\label{introduzione}

Toric quasifolds are highly singular spaces that extend simplicial toric varieties to nonrational fans.
They were first introduced in the symplectic category \cite{pcras,p}, by 
means of a generalization of the classical Delzant
construction of symplectic toric manifolds \cite{delzant}. 
Complex toric quasifolds, on the other hand,
first appeared in \cite{cx} and follow the construction of 
simplicial toric varieties as complex quotients given by Audin in
\cite{audin}. Toric quasifolds provide an interesting class of examples of 
diffeological spaces: the notion of diffeological quasifold was introduced, jointly with
Iglesias--Zemmour, in \cite{IZP} and
related to that of quasifold groupoid by Karshon--Miyamoto in \cite{KM}.

Quasifolds are locally the quotient of a manifold modulo the action of a countable group. 
If the countable groups are all finite, quasifolds are orbifolds; if they are all trivial, quasifolds are manifolds. 
Similarly to what happens for manifolds, in order to define a quasifold structure 
we need an atlas, namely a cover by mutually compatible local models.
We have shown in \cite{cx} that every $n$--dimensional complex toric quasifold
has a finite atlas whose local models are given 
by $\C^n$ modulo the action of a countable subgroup of $(S^1)^n$. 
Now, recall from Danilov \cite[Section 0.2]{danilov} that a smooth toric variety is 
characterized by the fact that it is covered by a finite number of copies of $\C^n$, 
with changes of charts expressed in terms of Laurent monomials. 
In this article, we introduce generalized Laurent monomials and we prove that
the changes of charts of toric quasifolds can be expressed, similarly to what happens in the
smooth case, in terms of generalized Laurent monomials.

This allows for the possibility of defining toric quasifolds directly, by starting with a 
collection of complex quasifold models and by gluing them together by means of 
generalized Laurent monomials, thus extending what was originally done for smooth toric varieties
\cite{demazure}. We plan to address this matter and its implications in future work.

The article is structured as follows. In Section~\ref{laurent}, we define generalized Laurent monomials for
local quasifold models. In Section~\ref{atlas},
we recall the construction of a complex toric quasifold, we describe the canonical affine atlas 
and write the chart transitions in terms of generalized Laurent monomials. 
Finally, in Section~\ref{examples}, we compute
the monomials in a number of interesting examples in dimensions $1$, $2$, and $3$.

\section{Generalized Laurent monomials}\label{laurent}

The goal of this section is to show that, in the quasifold setting, we can make sense of the notion of Laurent 
monomials with real exponents. 

Let us consider, first of all, an example of countable group that will be of particular relevance for us,
that of a {\em quasilattice}. We recall from \cite{pcras, p} that a quasilattice $Q$ in $\R^n$
is an additive subgroup of $\R^n$ given by the $\Z$--span of a finite generating subset of $\R^n$. 
Notice that a quasilattice is a lattice if, and only if, it is generated by a basis of $\R^n$. 

Let us now describe a quasifold of one chart which
provides the local model for toric quasifolds. First of all, 
write the exponential mapping from the Lie algebra $\C^m$ to the torus $(\C^*)^m$ as follows: 
$$\begin{array}{ccccc}
\exp\colon&\C^m&\longrightarrow&(\C^*)^m\\
&(z_1,\ldots,z_m)&\longmapsto&(e^{2\pi i z_1},\ldots,e^{2\pi i z_m}).
\end{array}
$$
We denote $\exp$ also its restriction to the Lie algebra $\R^m$ of $\R^m/\Z^m=(S^1)^m$.

If $Q\subset\R^n$ is a quasilattice, then $\Gamma=\exp(Q)$ is a countable subgroup of $\R^n/\Z^n=(S^1)^n$. 
We consider the quasifold on one chart
$$\C^n/\Gamma.$$

We are now ready to introduce the generalized Laurent monomials.
Consider $\underline{a}=(a_1,\ldots, a_n)\in Q$. If $l\in\{1,\ldots,n\}$
and $[z_1:\cdots:z_n]\in(\C^*)^n/\Gamma$,
then the class
$[z_l^{a_1}:\cdots:z_l^{a_n}]$
is well defined in $(\C^*)^n/\Gamma$. In fact, if we
write 
$$z_l=e^{2\pi i (x_l+h_l+iy_l)},$$
with $h_l\in\Z$,
then
$$z_l^{a_i}=e^{2\pi i a_i(x_l+h_l+iy_l)}.$$
Notice now that $\exp(h_l\underline{a})\in\Gamma$.

Similarly, consider a subset of indices $I\subseteq \{1,\ldots,n\}$ of cardinality $k$
and take $\underline{a}^l=(a_{1l},\ldots,a_{nl})\in Q$, $l\in I$. Then the class
$$\left[\prod_{l\in I}z_l^{a_{1l}}:\cdots:\prod_{l\in I}z_l^{a_{nl}}\right]$$
 is well defined in $(\C^*)^n/\Gamma$. 
 We call the single terms
 $$\prod_{l\in I} z_l^{a_{jl}}, \quad j=1,\ldots,n$$
 {\em generalized Laurent monomials}. 

\section{The canonical atlas and the generalized Laurent monomials}
\label{atlas}
\subsection{Fans}
Let us recall some basic facts on fans. For further details, we refer the reader to Cox et al. \cite{cox}
and Ziegler \cite{ziegler}.

A fan $\Sigma$ in $\R^n$ is a collection of strongly convex 
polyhedral cones in $\R^n$, such that each nonempty face of a cone in $\Sigma$ is itself a cone in $\Sigma$, 
and such that the intersection of any two cones in $\Sigma$ is a face of each.

The cones of dimension $n$ are called \textit{maximal cones}, 
while the cones of dimension $1$ are called \textit{generating rays}. 

An important example of fan
is the \textit{normal fan} $\Sigma_{\Delta}$ of a convex polytope $\Delta$. 
It is the fan whose generating rays are inward pointing and orthogonal to the polytope facets. 
The rays are grouped to generate all the cones of the fan following the combinatorics 
of the polytope faces and we have an inclusion--reversing bijection between 
cones in $\Sigma_{\Delta}$ and faces of $\Delta$.

Each cone of a fan is spanned by a subset of a set of generating vectors
$X_1,\ldots,X_d$. The fan is called \textit{simplicial} if, 
for each cone, these vectors are linearly independent. The normal fan of a convex polytope 
$\Sigma_{\Delta}$ is simplicial if, and only if, the polytope $\Delta$ is simple.
Consider a simplicial fan. Then, for each maximal cone $\sigma$, there is a subset 
$I_\sigma\subset\{1,\ldots,d\}$ of $n$ indices, such that this maximal cone is generated by $X_i,i\in I_\sigma$.  
The faces of $\sigma$ are in bijective correspondence with the subsets of the set of indices 
$I_\sigma$, the $0$--dimensional cone of $\Sigma$ corresponding to the empty set.

The support $|\Sigma|$ of a fan $\Sigma$ is the union in $\R^n$ of its cones. We will be interested in fans
that have convex support of full dimension; this implies that $|\Sigma|$ is the union of the maximal cones of $\Sigma$. 
A notable example is the case of a \textit{complete} fan, 
namely a fan whose support is the whole space $\R^n$ like, for instance, the normal fan 
of a convex polytope. 

A fan $\Sigma$ in $\R^n$ is said to be \textit{rational} if there exists a lattice 
$L\subset\R^n$ such that the intersection of $L$ with each generating ray of $\Sigma$ is nonempty. 
A convex polytope $\Delta$ is said to be rational if its normal fan $\Sigma_{\Delta}$ is.

We are interested in fans and polytopes that are not rational. In this context, it is convenient
to introduce a different notion. A fan $\Sigma$ in $\R^n$ is said to be \textit{quasirational}
with respect to a quasilattice $Q\subset \R^n$ if the intersection of $Q$ with each 
generating ray of $\Sigma$ is nonempty. 
A convex polytope $\Delta$ is said to be quasirational if its normal fan $\Sigma_{\Delta}$ is.

\subsection{The construction of complex toric quasifolds}
We recall from \cite[Theorem 2.2]{cx} the construction of toric quasifolds in the complex category.
In the original paper, we had considered simple convex polytopes, but the same argument goes through
in the more general setting of \textit{simplicial fans with convex support of full dimension}.

Let $Q$ be a quasilattice in $\R^n$, 
and assume that $\Sigma$ is quasirational with respect to $Q$. Let $X_1,\ldots,X_d\in \R^n$ 
be a choice of generators for the $d$ rays of $\Sigma$ and assume that they lie in $Q$. 
In the case of the normal fan of a convex polytope $\Delta$, these vectors 
are each orthogonal to a different codimension one face (\textit{facet}) of $\Delta$ and point inwards;
we will be referring to them as \textit{normals} for $\Delta$.

The \textit{fundamental triple} 
$$(\Sigma, Q, \{X_1,\ldots,X_d\})$$
effectively replaces the standard triple
$$
\mbox{(fan, lattice, primitive generators)}
$$
that is needed for the construction of toric varieties in the rational case. 
We will in fact construct the $n$--dimensional complex toric quasifold, 
$X_\Sigma$, relatively to the fundamental triple above.

Before we proceed, we remark that, in this framework, tori are generalized by {\em quasitori}
\cite{pcras,p,cx}. They are given by the quasifolds of one chart
$$\R^n/Q$$ and $$\C^n/Q.$$
The first is real, while the second is complex. They are both abelian groups and, when 
$Q$ is a true lattice, $\R^n/Q$ and $\C^n/Q$ are true tori. Their Lie algebras are
given by $\R^n$ and $\C^n$, respectively.

Let now $\sigma$ be a maximal cone, and let $I_\sigma$ denote the corresponding index set.
Consider the open subset of $\C^d$ given by
$$U_{\sigma}=\{\,\underline{z}\in\C^d\mid z_j\neq0,\;\;j\notin I_\sigma\,\}.$$
Notice that we have
$$U_{\sigma}=\bigsqcup_{J\subset I_{\sigma}}\{\,\underline{z}\in\C^d\mid z_j=0\;\;\text{iff}\;\;j\in J\,\}.$$
Take now
$$U_{\Sigma}=\bigcup_{\sigma\;\text{maximal cone}}U_{\sigma}.$$
This open subset of $\C^d$ is acted upon by a subgroup $N_{\C}\subset (\C^*)^d$ 
that we construct in the following way.

Let $\{e_1,\ldots,e_d\}$ denote the standard basis
of $\R^d$ and $\C^d$; consider the surjective linear mapping
$$
\begin{array}{cccc}\label{pi}
\pi \,\colon\,& \R^d & \longrightarrow & \R^n\\
&   e_j& \longmapsto & X_j
\end{array}
$$ 
and its complexification $$
\begin{array}{cccc}
\pi_{\C} \,\colon\,& \C^d & \longrightarrow & \C^n\\
&   e_j& \longmapsto & X_j.
\end{array}
$$ Consider the quasitorus $\R^n/Q$ and its complexification
$\C^n/Q$. The mappings $\pi$ and $\pi_{\C}$ induce group homomorphisms
$$\Pi
\,\colon\, (S^1)^d=\R^d/\Z^d\longrightarrow \R^n/Q$$ and $$\Pi_{\C} \,\colon\, (\C^*)^d=\C^d/\Z^d\longrightarrow \C^d/Q.$$
We define $N$ to be the kernel of the mapping $\Pi$ and $N_{\C}$ to be the kernel of the
mapping $\Pi_{\C}$. Notice that neither $N$ nor $N_{\C}$ are honest tori unless $Q$ is a honest
lattice.

Then the complex toric quasifold corresponding the fundamental triple $(\Sigma,Q,\{X_1,\ldots,X_d\})$ is
given by the quotient
$$X_\Sigma=\frac{U_\Sigma}{N_{\C}}.$$
It is an $n$--dimensional complex quasifold.
If $\Sigma$ is the normal fan of a convex polytope $\Delta$, we will call this quasifold 
$X_\Delta$.

Notice that the mapping $\Pi_{\C}$ induces an isomorphism
$$
(\C^*)^d/N_{\C}\longrightarrow \C^n/Q.
$$
Thus the complex quasitorus $\C^n/Q$ acts on $X_\Sigma$. This action is holomorphic and, as in the rational
case, $X_\Sigma$ is the disjoint union of the $k$--dimensional orbits, each corresponding to an 
$(n-k)$--dimensional cone of $\Sigma$. In particular, there is a dense open orbit, 
corresponding to the $0$--dimensional cone.
Also, there are a finite number of fixed points, one for each maximal cone $\sigma$. They are given by
the points $[z_1:\cdots:z_d]$, with $z_j=0$, for $j\in I_\sigma$ and $z_j=1$ otherwise.

\begin{rem}\rm{
In recent years, nonrational toric geometry has been investigated from several different, though
interrelated, points of view. See \cite{whatis} for a detailed account.}
\end{rem}

\subsection{The symplectic picture}
It was shown in \cite{pcras, p} that, to each fundamental triple
$$(\Sigma_\Delta, Q, \{X_1,\ldots,X_d\}),$$ with $\Delta$ a simple convex polytope, there corresponds a 
$2n$--dimensional symplectic toric quasifold $M_\Delta$, endowed with the effective Hamiltonian action of the
real quasitorus $\R^n/Q$. The image of the corresponding moment mapping $\Phi$ is the initial polytope
$\Delta$. The quasifold $M_\Delta$ is obtained explicitly
 via symplectic reduction with respect to the subgroup
$N\subset (S^1)^d$ that we have introduced above, thereby extending the classical
construction of symplectic toric manifolds by Delzant \cite{delzant}.
It is shown in \cite[Theorem 3.2]{cx} that $M_\Delta$ and 
$X_\Delta$ are equivariantly diffeomorphic, and that the complex and symplectic structures are compatible, 
thus defining a Kahler structure on $X_\Delta$, exactly as in the rational case (for which we refer to Audin \cite{audin}
and Guillemin \cite{guillemin}). 
Moreover, the fixed points 
of the action of $\C^n/Q$ on $X_\Delta$
correspond to the fixed points of the action of $\R^n/Q$ on $M_\Delta$, and the latter are mapped bijectively via $\Phi$ to
the vertices of the polytope.

\subsection{The canonical affine atlas and the generalized Laurent monomials}
We adapt from \cite{cx} the construction of the canonical atlas for $X_\Sigma$. 
The novelty here, besides the weaker assumptions on $\Sigma$, 
is that we write the changes of charts explicitly in terms of 
generalized Laurent monomials.

Before we do so, notice that we can write the groups $N$ and $N_{\C}$ as follows:
$$N=\{\,\exp(X)\mid X\in \R^d,\;\;\pi(X)\in Q\,\}$$
$$N_{\C}=\{\,\exp(Z)\mid Z\in \C^d,\;\;\pi(Z)\in Q\,\}.$$
The respective Lie algebras are
$\mathfrak{n}=\ker(\pi)$ and
$\mathfrak{n}_{\C}=\ker(\pi_{\C})$, which are linear subspaces of dimension $d-n$ of $\R^d$ and $\C^d$, respectively. 
The mapping $\exp$ restricts to both.

Fix a maximal cone, $\sigma$, with its subset of $n$ indices, $I_{\sigma}$.
Consider the restrictions $\pi_{\sigma}$ and $(\pi_{\sigma})_\C$ of the mappings $\pi$ and
$\pi_{\C}$ to the $n$--dimensional subspaces $\prod_{i\in I_\sigma}\R e_i\subset\R^d$ and 
$\prod_{i\in I_\sigma}\C e_i\subset \C^d$ respectively. Both maps are isomorphisms.
We have the following lemma, that, in the complex setting, is proven in \cite[Lemma 2.3]{cx}. 
We give a new proof that is tailored to our purposes:
\begin{lem} For each maximal cone $\sigma$, we have that $N=\Gamma_\sigma \exp(\mathfrak{n})$
and $N_{\C}=\Gamma_\sigma \exp(\mathfrak{n}_\C)$, where 
$$\Gamma_\sigma=\left\{\,\exp(X)\mid X\in\prod_{j\in I_\sigma}\R e_j,\;\;\pi(X)\in Q\,\right\}.$$
\end{lem}
\noindent\begin{preuve}
Let $X\in\R^d$ such that $\pi(X)\in Q$ and let $Y=\pi_\sigma^{-1}(\pi(X))$. Then $W=X-Y\in \mathfrak{n}$ 
and therefore $\exp(X)=\exp(Y)\exp(W)$, with $\exp(Y)\in \Gamma_\sigma$,
$\exp(W)\in \exp(\mathfrak{n})$. The same argument applies to the 
complexified group and Lie algebra.
\end{preuve}

\begin{rem}{\rm
For each $j\notin I_{\sigma}$ denote $\underline{a}^j=\pi_{\sigma}^{-1}(X_j)\in \prod_{i\in I_\sigma}\R e_i\subset\R^d$. The $d-n$ linearly independent vectors $e_j-\underline{a}^j$ form a basis of $\mathfrak{n}$ and $\mathfrak{n}_{\C}$.
Notice that $\exp(\underline{a}^j)\in \Gamma_{\sigma}$.}
\end{rem}
We are finally ready to prove the following
 \begin{thm} The complex toric quasifold $X_\Sigma$ corresponding to the fundamental triple 
 $(\Sigma,Q,\{X_1,\ldots,X_d\})$ is endowed with a generalized canonical atlas, having
a chart around each of the fixed points of the quasitorus action. 
The transition mappings are written in terms of generalized Laurent monomials.
\end{thm} 
\noindent\begin{preuve}
We define a chart for each maximal cone in the fan $\Sigma$. Given such a cone,
$\sigma$, we can assume, up to renumbering, that the corresponding index set, $I_\sigma$,
is given by $\{1,\ldots,n\}$. Then the mapping
$$
\begin{array}{cccc}
&\C^n/\Gamma_\sigma & \stackrel{\eta_\sigma}\longrightarrow & V_\sigma=U_\sigma/N_{\C} \\
&[z_1: \cdots : z_n] & \longmapsto & [z_1:\cdots:z_n:1:\cdots:1]
\end{array}
$$
is a homeomorhism and defines a chart around the fixed point 
$$[0:\cdots:0:1:\cdots:1].$$
We will be referring to this chart as the quadruple $(\C^n,\Gamma_\sigma,V_\sigma,\eta_\sigma)$.

Consider now a second chart $(\C^n,\Gamma_\tau,V_\tau,\eta_\tau)$, 
corresponding to another maximal cone $\tau$, and assume that $V_\sigma\cap V_\tau\neq\emptyset$. 
This happens whenever $I_\sigma\cap I_\tau\neq\emptyset$. 
Assume, again up to renumbering, that we have the following partition 
$$\{1,\ldots,h,h+1,\ldots,n,n+1,\ldots,n+h,n+h+1,\ldots,d\}, \;\;1\leq h \leq n-1,$$ with
$$I_{\sigma}=\{1,\ldots,h,h+1,\ldots,n\},$$
$$I_{\tau}=\{h+1,\ldots,n,n+1,\ldots,n+h\},$$
$$I_{\sigma}\cap I_{\tau}=\{h+1,\ldots,n\}.$$
Let us compute the transition mapping between these two charts.
Take
$$[z_{h+1}:\cdots:z_{n+h}]\stackrel{\eta_\tau}\longmapsto[1:\cdots:1:z_{h+1}:\cdots:z_{n+h}:1:\cdots:1].$$
Let us now move the representative on the right--hand side within the $N_{\C}$--orbit by acting via
$$
\left(\prod_{j=n+1}^{n+h}z_{j}^{a_{1j}},\ldots,\prod_{j=n+1}^{n+h}z_{j}^{a_{nj}},
z_{n+1}^{-1},\ldots,z_{n+h}^{-1},1,\ldots,1\right)\in\exp(\mathfrak{n}_{\C}),
$$
where $\underline{a}^j=\pi_{\sigma}^{-1}(X_j)=(a_{1j},\ldots,a_{nj})$.
We obtain
$$\left[\prod_{j=n+1}^{n+h}z_{j}^{a_{1j}}:\cdots:\prod_{j=n+1}^{n+h}z_{j}^{a_{hj}}:
	\prod_{j=n+1}^{n+h}z_{j}^{a_{h+1,j}}z_{h+1}:\cdots:\prod_{j=n+1}^{n+h}z_{j}^{a_{nj}}z_n:
	1:\cdots:1\right].$$
Write now
$\prod_{j=n+1}^{n+h}z_{j}^{a_{ij}}=\underline{z}^{\underline{a}_i},$
where $\underline{z}=(z_{n+1},\ldots,z_{n+h})\in\C^h$ and
$$\underline{a}_{i}=(a_{i,n+1},\ldots,a_{i,n+h})\in\R^h,\;\;
i=1,\ldots,n.$$

Then the transition mapping from $(\C^n,\Gamma_\tau,V_\tau,\eta_\tau)$ to $(\C^n,\Gamma_\sigma,V_\sigma,\eta_\sigma)$ can be finally written in terms of generalized Laurent monomials as
$$
\begin{array}{ccc}
\eta_\tau^{-1}(V_\sigma\cap V_\tau)&\stackrel{\eta_\sigma^{-1}\circ\eta_\tau}\longrightarrow&\eta_\sigma^{-1}(V_\sigma\cap V_\tau)\\

[z_{h+1}:\cdots:z_{n+h}]&\longmapsto&[\underline{z}^{\underline{a}_1}:\cdots:\underline{z}^{\underline{a}_h}:\underline{z}^{\underline{a}_{h+1}}z_{h+1}:\cdots:\underline{z}^{\underline{a}_{n}}z_n].
\end{array}
$$
\end{preuve}

\section{Examples}\label{examples}
We compute the generalized Laurent monomials in a number of interesting examples.
It will be convenient to label the charts relatively to the index sets $I_\sigma$ instead 
of the maximal cones $\sigma$ themselves.
\subsection{The quasisphere}\label{quasisfera}
The quasisphere is a toric quasifold that 
was first introduced in \cite{pcras,p} and can be constructed as follows.
Let $a$ be a positive irrational number and consider the unit interval $[0,1]$.
Then the quasisphere is the toric quasifold obtained from the fundamental triple
$$
(\Sigma_{[0,1]}, \Z+a\Z, \{X_1=a, X_2=-1\}).
$$
It is given by
$$
X_a=\frac{\mathbb{C}^2\setminus\{0\}}{\{\,(e^{2\pi i r}, e^{2\pi i ar})\in(\mathbb{C}^*)^2\mid r\in \mathbb{C}\,\}}
$$
and it is endowed with the holomorphic action of the complex dimension one quasitorus (quasicircle)
$$\mathbb{C}/(\mathbb{Z}+a\mathbb{Z}).$$
Notice that this quasicircle is the complexification of the {\em irrational torus} of Donato--Iglesias \cite{di}.
The canonical atlas of the complex quasisphere is made of two charts, one around the point
$[0:1]$, the other around the point $[1:0]$.
The first is given by $(\mathbb{C},\Gamma_1, V_1, \eta_1)$, where
$$
\Gamma_1=\left\{\,e^{2\pi i \frac{h}{a}}\in S^1\mid h\in\mathbb{Z}\right\},
$$
$$V_1=\{[z_1:z_2] \in X_a\mid z_2\neq0\},$$
$$
\begin{array}{ccc}
\mathbb{C}/\Gamma_1&\stackrel{\eta_1}{\longrightarrow}& V_1\\
\,[z]&\longmapsto&\left[ z:1\right].
\end{array}
$$
While the second
is given by $(\mathbb{C},\Gamma_2, V_2, \eta_2)$, where 
$$
\Gamma_2=\left\{\,e^{2\pi i ah}\in S^1\mid h\in\mathbb{Z}\right\},
$$
$$V_2=\{[w_1:w_2] \in X_a\mid w_1\neq0\},$$
$$
\begin{array}{ccc}
\C/\Gamma_2&\stackrel{\eta_2}{\longrightarrow}& V_2\\
\,[w]&\longmapsto&\left[1:w\right].
\end{array}
$$
The transition map between these two charts is given by
$$
\begin{array}{ccc}
\eta_1^{-1}(V_1\cap V_2)&\stackrel{\eta_2^{-1}\circ\eta_1}\longrightarrow& \eta_2^{-1}(V_1\cap V_2)\\
\,[z]&\longmapsto&\left[z^{-a}\right].
\end{array}
$$
Here we get the simplest possible Laurent monomial. 

A thorough description of the quasisphere from the symplectic viewpoint can be found in \cite{p4}.

\subsection{Generalized weighted projective space}
\begin{figure}[h]
 \begin{center}
\includegraphics[scale=0.6]{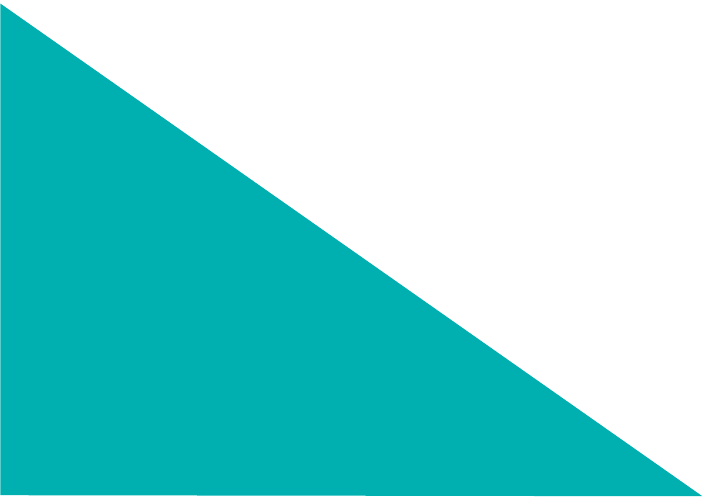} \quad\quad
\includegraphics[scale=0.2]{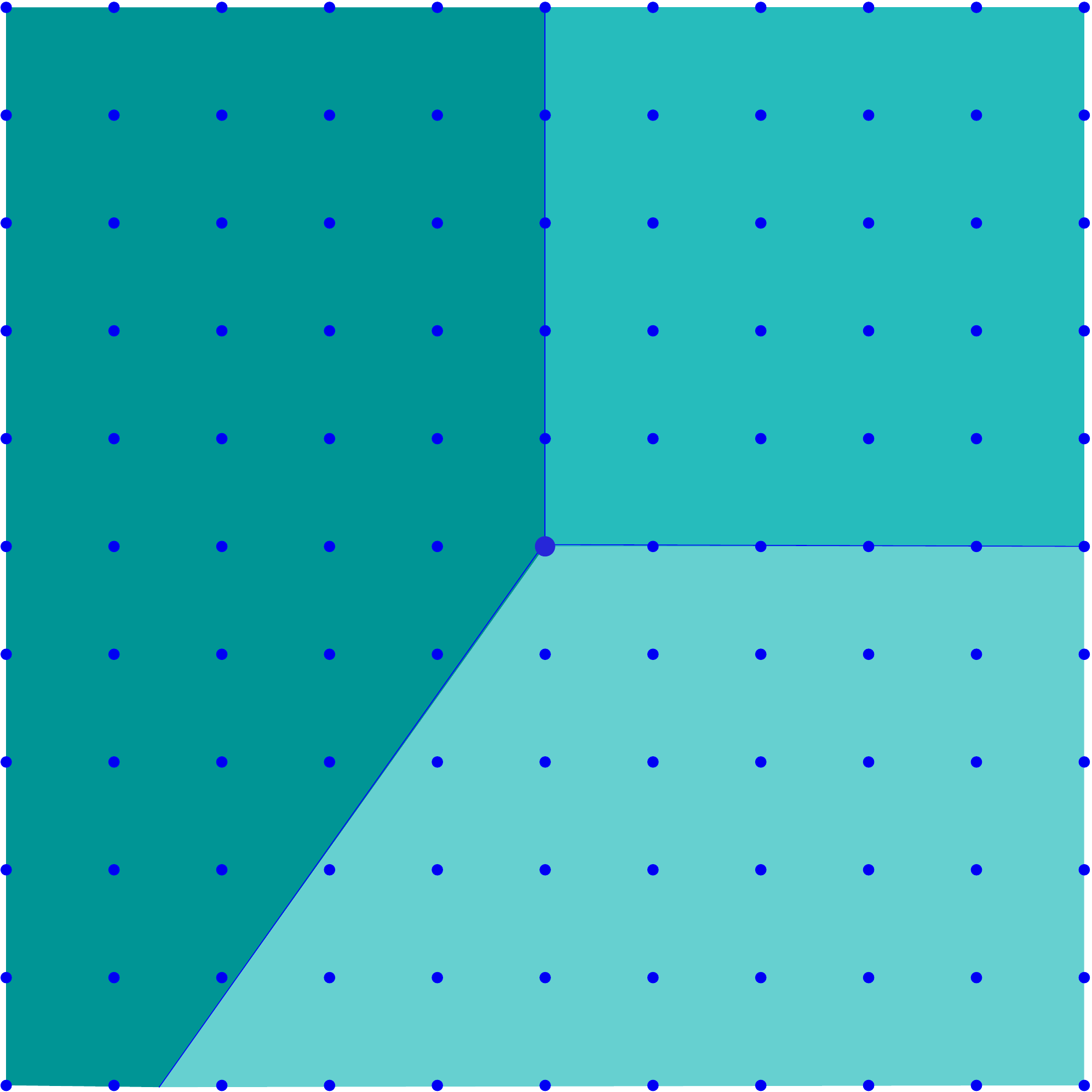}
\caption{The right triangle $T_a$ and its normal fan.}
\label{triangle}
\end{center}
\end{figure}

Consider, for any positive real number $a$, the right triangle $T_a$ of vertices
$(0,0)$, $(a,0)$, and $(0,1)$ (see Figure~\ref{triangle}). 

Though this triangle is rational with respect to the lattice spanned by
$(1,0)$ and $(0,a)$, it is interesting to view it in a nonrational setting, as follows. Let us take
normals given by 
$X_1=(1,0)$, $X_2=-(1,a)$ and $X_3=(0,1)$. They span the quasilattice
$$\mathbb{Z}\times (\mathbb{Z}+a\mathbb{Z})\supseteq \mathbb{Z}^2.$$ 
The toric quasifold corresponding to the fundamental triple
$$
(\Sigma_{T_a}, \mathbb{Z}\times (\mathbb{Z}+a\mathbb{Z}), \{X_1, X_2, X_3\})
$$
is given by
$$\C\P^2_{(1,1,a)}=\frac{\mathbb{C}^3\setminus\{0\}}{\{\,(e^{2\pi i u}, e^{2\pi i u}, e^{2\pi i au})\in(\mathbb{C}^*)^3\mid u\in \mathbb{C}\,\}}.$$
It is endowed with the holomorphic action of the complex quasitorus
$$\mathbb{C}^2/[\mathbb{Z}\times (\mathbb{Z}+a\mathbb{Z})]\simeq \mathbb{C}^*\times (\mathbb{C}/(\mathbb{Z}+a\mathbb{Z})).$$
Notice that for $a=1$ we get ordinary complex projective space $\C\P^2$, while for $a=n$, a positive integer, we get 
the weighted projective space $\C\P^2_{(1,1,n)}$. 
This complex toric quasifold admits an atlas with three charts. Consider the charts around the fixed points $[1:0:0]$
and $[0:1:0]$. The first is given
 by $(\mathbb{C}^2,\Gamma_{23}, V_{23}, \eta_{23})$, where 
$$
\Gamma_{23}=\left\{\,(1,e^{2\pi ia h})\in (S^1)^2\mid h\in\mathbb{Z}\right\},
$$
$$V_{23}=\{[z_1:z_2:z_3] \in \C\P^2_{(1,1,a)}\mid z_1\neq 0\},$$
$$
\begin{array}{ccc}
\mathbb{C}^2/\Gamma_{23}&\stackrel{\eta_{23}}{\longrightarrow}& V_{23}\\
\,[z_2:z_3]&\longmapsto&\left[ 1:z_2:z_3\right].
\end{array}
$$
The second is given by
is given by $(\mathbb{C}^2,\Gamma_{13}, V_{13}, \eta_{13})$, where 
$$
\Gamma_{13}=\left\{\,(1,e^{2\pi ia h})\in (S^1)^2\mid h\in\mathbb{Z}\right\},
$$
$$V_{13}=\{[w_1:w_2:w_3] \in \C\P^2_{(1,1,a)}\mid w_2\neq0\},$$
$$
\begin{array}{ccc}
\mathbb{C}^2/\Gamma_{13}&\stackrel{\eta_{13}}{\longrightarrow}& V_{13}\\
\,[w_1:w_3]&\longmapsto&\left[w_1:1:w_3\right].
\end{array}
$$
The transition map and Laurent monomials between these two charts are given by
$$
\begin{array}{ccc}
\eta_{23}^{-1}(V_{23}\cap V_{13})&\stackrel{\eta_{13}^{-1}\circ\eta_{23}}\longrightarrow& \eta_{13}^{-1}(V_{23}\cap V_{13})\\
\,[z_2:z_3]&\longmapsto&\left[z_2^{-1} : z_2^{-a}z_3\right].
\end{array}
$$

\subsection{Generalized Hirzebruch surfaces}
We recall the construction of a one--parameter family of toric quasifolds generalizing Hirzebruch surfaces which
was first introduced jointly with Zaffran in \cite{hirze}.

\begin{figure}[h]
\begin{center}
\includegraphics[scale=0.5]{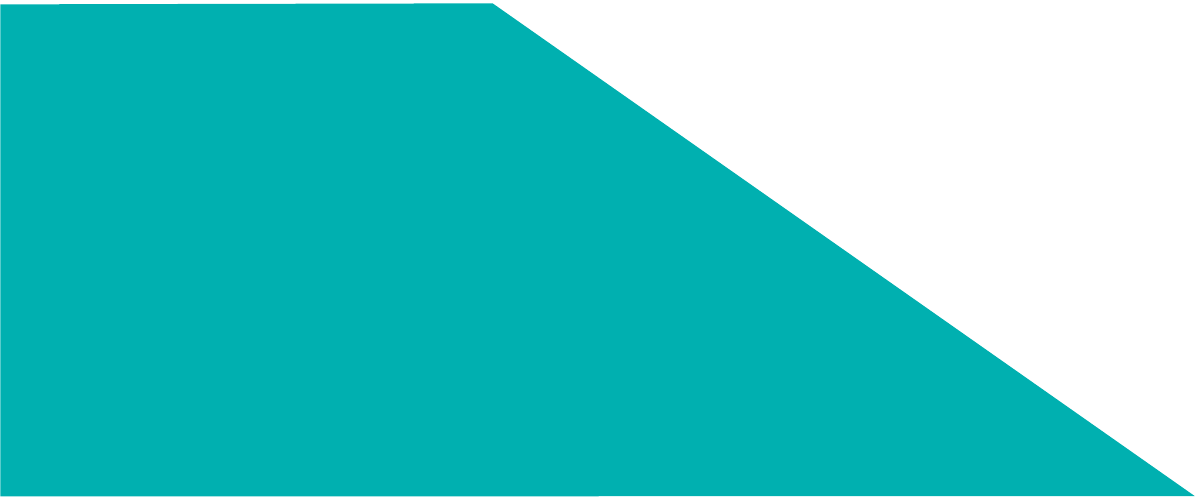} \quad\quad
\includegraphics[scale=0.2]{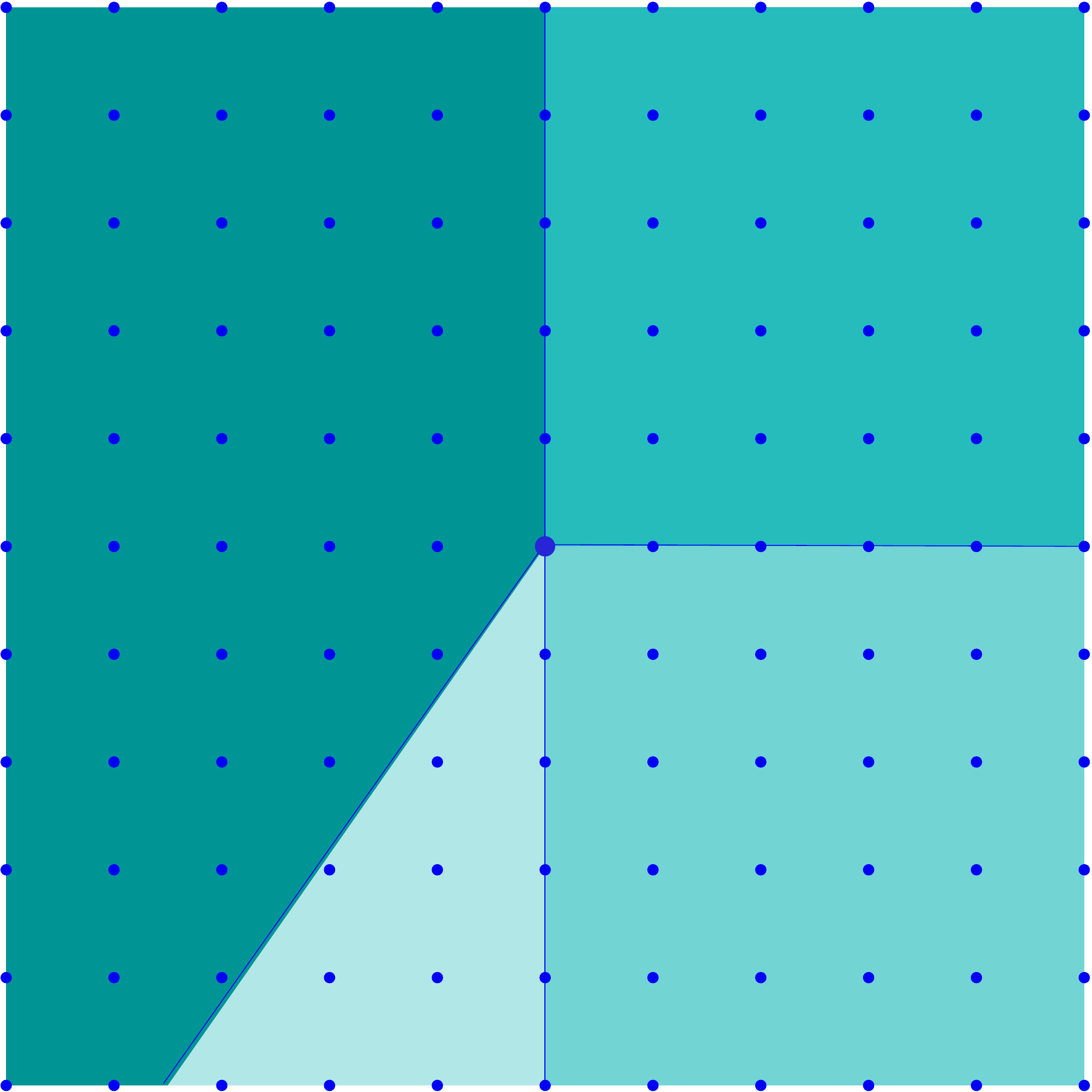}
\caption{The trapezoid $\mathcal{T}_a$ and its normal fan.}
\label{trapezoid}
\end{center}
\end{figure}

Consider, for any positive real number $a$, the trapezoid $\mathcal{T}_a$
of vertices $(0,0)$, $(1,0)$, $(1,1)$, and $(a+1,0)$ (see Figure~\ref{trapezoid}). 
We recall that $\mathcal{T}_n$, for $n$ a positive integer, 
is the polytope that corresponds to the Hirzebruch surface $\H_n$.
Similarly to what happens with the right triangle $T_a$, the trapezoid 
$\mathcal{T}_a$ is always rational with respect to the lattice generated by $(1,0)$, $(0,a)$.
However, again as above, we view $\mathcal{T}_a$ in a nonrational setting. To do so, we choose
normals given by $X_1=(1,0)$, $X_2=-(1,a)$, $X_3=(0,1)$, and $X_4=(0,-1)$;
they too span the quasilattice
$\mathbb{Z}\times (\mathbb{Z}+a\mathbb{Z})$.
The toric quasifold corresponding to the fundamental triple
$$
(\Sigma_{\mathcal{T}_a}, \mathbb{Z}\times (\mathbb{Z}+a\mathbb{Z}), \{X_1, X_2, X_3,X_4\}).
$$
is given, if $\underline{z}=(z_1,z_2,z_3,z_4)\in\mathbb{C}^4$, by
$$\H_a=\frac{\{\, \underline{z}\mid z_1z_3\neq0\,\}\cup\{\, \underline{z}\mid z_1z_4\neq0\,\}\cup\{\, \underline{z}\mid z_2z_3\neq0\,\}\cup\{\, \underline{z}\mid z_2z_4\neq0\,\}}{\{\,(e^{2\pi i u}, e^{2\pi i u}, e^{2\pi i (au+v)}, e^{2\pi i v})\in(\mathbb{C}^*)^4\mid u,v\in \mathbb{C}\,\}}=$$
$$
=\frac{(\mathbb{C}^2\setminus\{0\})\times(\mathbb{C}^2\setminus\{0\})}{\{\,(e^{2\pi i u}, e^{2\pi i u}, e^{2\pi i (au+v)}, e^{2\pi i v})\in(\mathbb{C}^*)^4\mid u,v\in \mathbb{C}\,\}}.
$$
It is also endowed with the holomorphic action of the complex quasitorus
$$\mathbb{C}^*\times (\mathbb{C}/(\mathbb{Z}+a\mathbb{Z})).$$
Notice that for $a=n$ we obtain the Hirzebruch surface $\H_n$.
This quasifold is covered by $4$ complex charts; 
let us describe the ones around the points $[1:0:0:1]$ and $[0:1:0:1]$.
The chart around the point $[1:0:0:1]$
is given by $(\mathbb{C}^2,\Gamma_{23}, V_{23}, \eta_{23})$, where 
$$
\Gamma_{23}=\left\{\,(1,e^{2\pi ia h})\in (S^1)^2\mid h\in\mathbb{Z}\right\},
$$
$$V_{23}=\{[z_1:z_2:z_3:z_4] \in \H_a\mid z_1z_4\neq 0\},$$
$$
\begin{array}{ccc}
\mathbb{C}^2/\Gamma_{23}&\stackrel{\eta_{23}}{\longrightarrow}& V_{23}\\
\,[z_2:z_3]&\longmapsto&\left[ 1:z_2:z_3:1\right].
\end{array}
$$
On the other hand, the chart around the point $[0:1:0:1]$
is given by $(\mathbb{C}^2,\Gamma_{13}, V_{13}, \eta_{13})$, where 
$$
\Gamma_{13}=\left\{\,(1,e^{2\pi ia h})\in (S^1)^2\mid h\in\mathbb{Z}\right\},
$$
$$V_{13}=\{[w_1:w_2:w_3:w_4] \in \H_a\mid w_2w_4\neq0\},$$
$$
\begin{array}{ccc}
\mathbb{C}^2/\Gamma_{13}&\stackrel{\eta_{13}}{\longrightarrow}& V_{13}\\
\,[w_1:w_3]&\longmapsto&\left[w_1:1:w_3:1\right].
\end{array}
$$
The transition map between these two charts and the corresponding Laurent monomials are given by
$$
\begin{array}{ccc}
\eta_{23}^{-1}(V_{23}\cap V_{13})&\stackrel{\eta_{13}^{-1}\circ\eta_{23}}
\longrightarrow& \eta_{13}^{-1}(V_{23}\cap V_{13})\\
\,[z_2:z_3]&\longmapsto&\left[z_2^{-1} : z_2^{-a}z_3\right].
\end{array}
$$
Notice that this change of charts is the same as the one for the weighted projective space above.
This is not a coincidence, but rather a manifestation of the fact that, analogously to what happens in the smooth case, 
the generalized Hirzebruch surface $\H_a$ is the blow--up of the generalized
weighted projective space $\C\P^2_{(1,1,a)}$ at its third fixed point $[0:0:1]$
(see \cite{cut, hirze}).

\subsection{The Penrose kite}
\begin{figure}[h]
 \begin{center}
\includegraphics[scale=0.4]{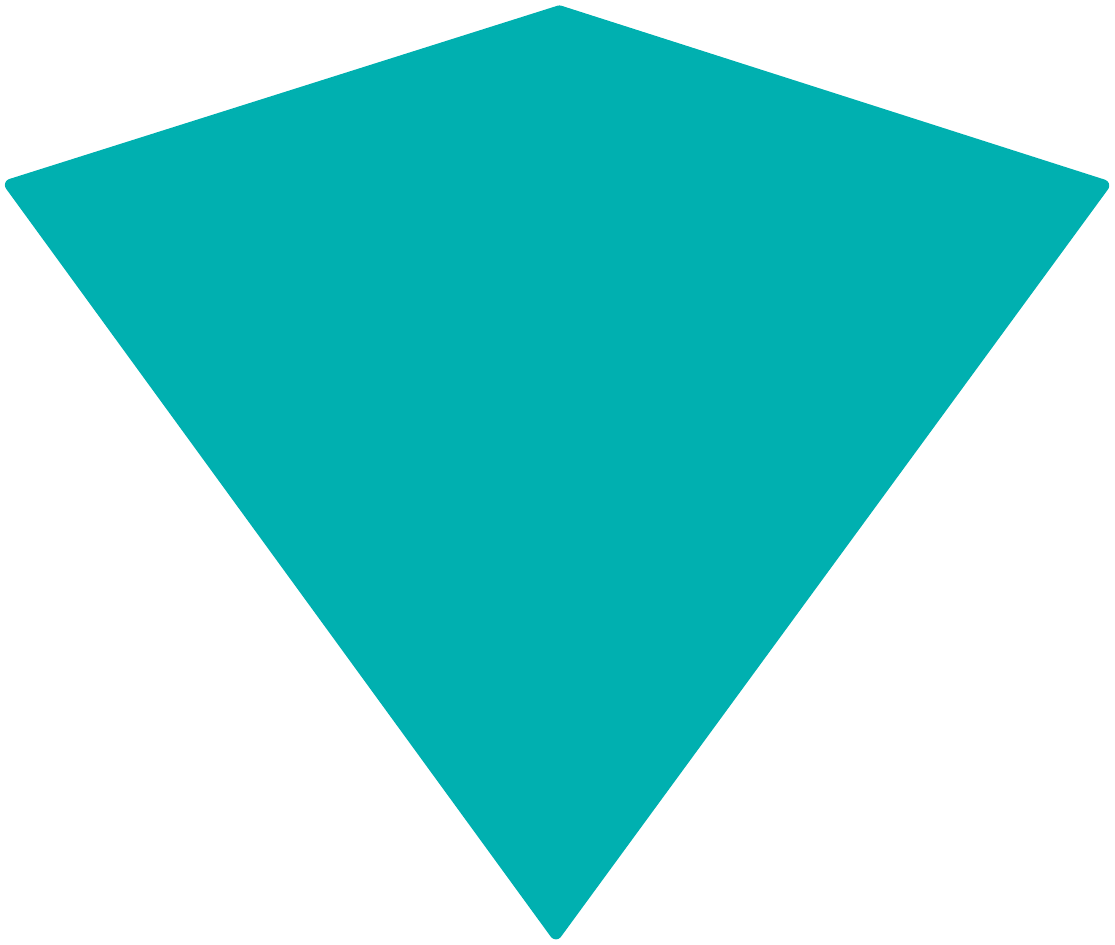} \quad\quad
\includegraphics[scale=0.2]{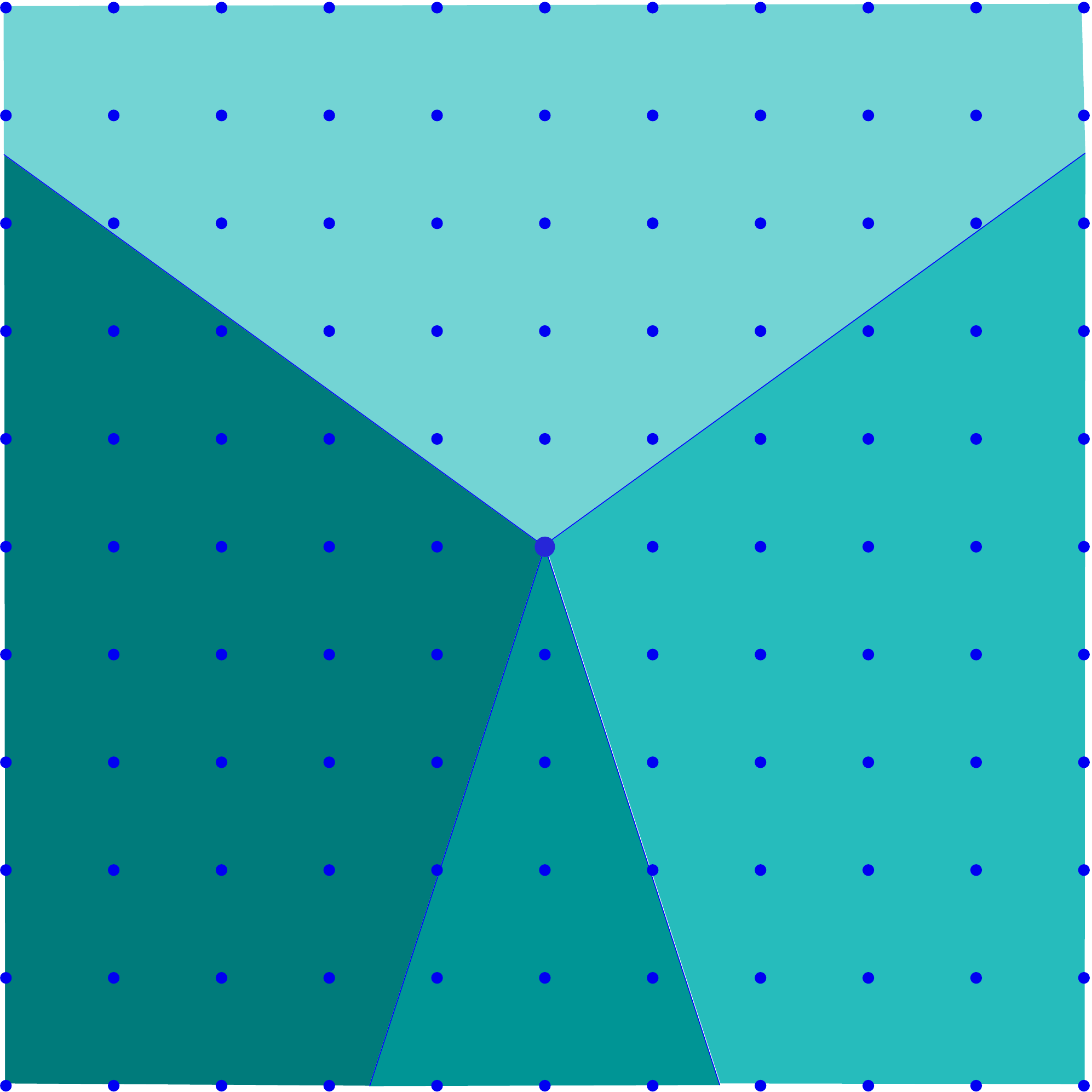}
\caption{The Penrose kite and its normal fan.}
 \label{kite}
 \end{center}
\end{figure}

The kite and dart tiling is an aperiodic tiling of the plane introduced by Penrose \cite{pentaplexity}. 
The kite is a quadrilateral having three of its angles equal to $\frac{2\pi}{5}$, the fourth equal to $\frac{4\pi}{5}$. Its
 long edges are $\phi$ times the short edges, where $\phi=\frac{1+\sqrt{5}}{2}$ is the
 {\em golden ratio} (see Figure~\ref{kite}).
 The kite is not rational, however it is quasirational with respect to the 
 pentagonal quasilattice $Q_5$ generated by the fifth--roots of unity:
$$
\begin{array}{l}
Y_0=(1,0)\\
Y_1=(\cos{\frac{2\pi}{5}},\sin{\frac{2\pi}{5}})=\frac{1}{2}(\frac{1}{\phi},\sqrt{2+\phi})\\
Y_2=(\cos{\frac{4\pi}{5}},\sin{\frac{4\pi}{5}})=\frac{1}{2}(-\phi,\frac{1}{\phi}\sqrt{2+\phi})\\
Y_3=(\cos{\frac{6\pi}{5}},\sin{\frac{6\pi}{5}})=\frac{1}{2}(-\phi,-\frac{1}{\phi}\sqrt{2+\phi})\\
Y_4=(\cos{\frac{8\pi}{5}},\sin{\frac{8\pi}{5}})=\frac{1}{2}(\frac{1}{\phi},-\sqrt{2+\phi}).
\end{array}
$$
Choose normals $X_1=-Y_1$, $X_2=-Y_3$, $X_3=Y_2$ and $X_4=Y_4$. Notice that
$$\phi=1+\frac{1}{\phi}$$ 
and that
$$\left\{
\begin{array}{lrrr}
Y_2=-&Y_1&-&\phi Y_4\\
Y_3=-&\phi Y_1&-& Y_4.\\
\end{array}
\right.
$$
The toric quasifold corresponding to the fundamental triple
$$
(\Sigma_{\Delta}, Q_5, \{X_1, X_2, X_3,X_4\})
$$
was first introduced and studied in \cite{kite}. It is given by
$$X_{\Delta}=\frac{(\mathbb{C}^2\setminus\{0\})\times(\mathbb{C}^2\setminus\{0\})}{\left\{\,\left(e^{2\pi i(\phi u-v)},e^{2\pi i u},e^{2\pi iv},e^{2\pi i(\phi v-u)}\right)\in (\C^*)^4\mid u,v \in\C\,\right\}}.$$
It is shown in \cite[Theorem 6.1]{kite} that {\em it is not} a global quasifold, 
namely not the quotient of a manifold modulo the action of a countable group. 

Again, here we have an atlas made of $4$ charts; let us describe the ones around the fixed points
$[0:1:1:0]$ and $[1:0:1:0]$.
The first is given by $(\C^2,\Gamma_{14}, V_{14}, \tau_{14})$, where
$$
\Gamma_{14}=\left\{\,(e^{2\pi i\phi h},e^{2\pi i\phi k})\in (S^1)^2\mid  h,k\in\mathbb{Z}\right\},
$$
$$V_{14}=\{[z_1:z_2:z_3:z_4] \in X_{\Delta}\mid z_2z_3\neq0\},$$
$$
\begin{array}{ccc}
\C^2/\Gamma_{14}&\stackrel{\tau_{14}}{\longrightarrow}& V_{14}\\
\,[z_1:z_4]&\longmapsto&\left[z_1:1: 1:z_4\right].
\end{array}
$$
The second is given by $(\C^2,\Gamma_{24}, V_{24}, \tau_{24})$, where
$$
\Gamma_{24}=\left\{\,(e^{2\pi i\phi h}e^{2\pi i\phi k},e^{-2\pi i\phi h})\in (S^1)^2\mid  h,k\in\mathbb{Z}\right\},
$$
$$V_{24}=\{[z_1:z_2:z_3:z_4] \in X_{\Delta}\mid z_1z_3\neq0\},$$
$$
\begin{array}{ccc}
\C^2/\Gamma_{24}&\stackrel{\tau_{24}}{\longrightarrow}& V_{24}\\
\,[w_2:w_4]&\longmapsto&\left[1:w_2:1:w_4\right].
\end{array}
$$
Their transition and Laurent monomials are given by
$$
\begin{array}{ccc}
\eta_{14}^{-1}(V_{14}\cap V_{24})&\stackrel{\eta_{24}^{-1}\circ\eta_{14}}
\longrightarrow& \eta_{24}^{-1}(V_{14}\cap V_{24})\\
\,[z_1:z_4]&\longmapsto&\left[z_1^{-\frac{1}{\phi}} : z_1^{\frac{1}{\phi}}z_4\right].
\end{array}
$$

\subsection{The regular dodecahedron}
\begin{figure}[h]
\begin{center}
\includegraphics[scale=0.25]{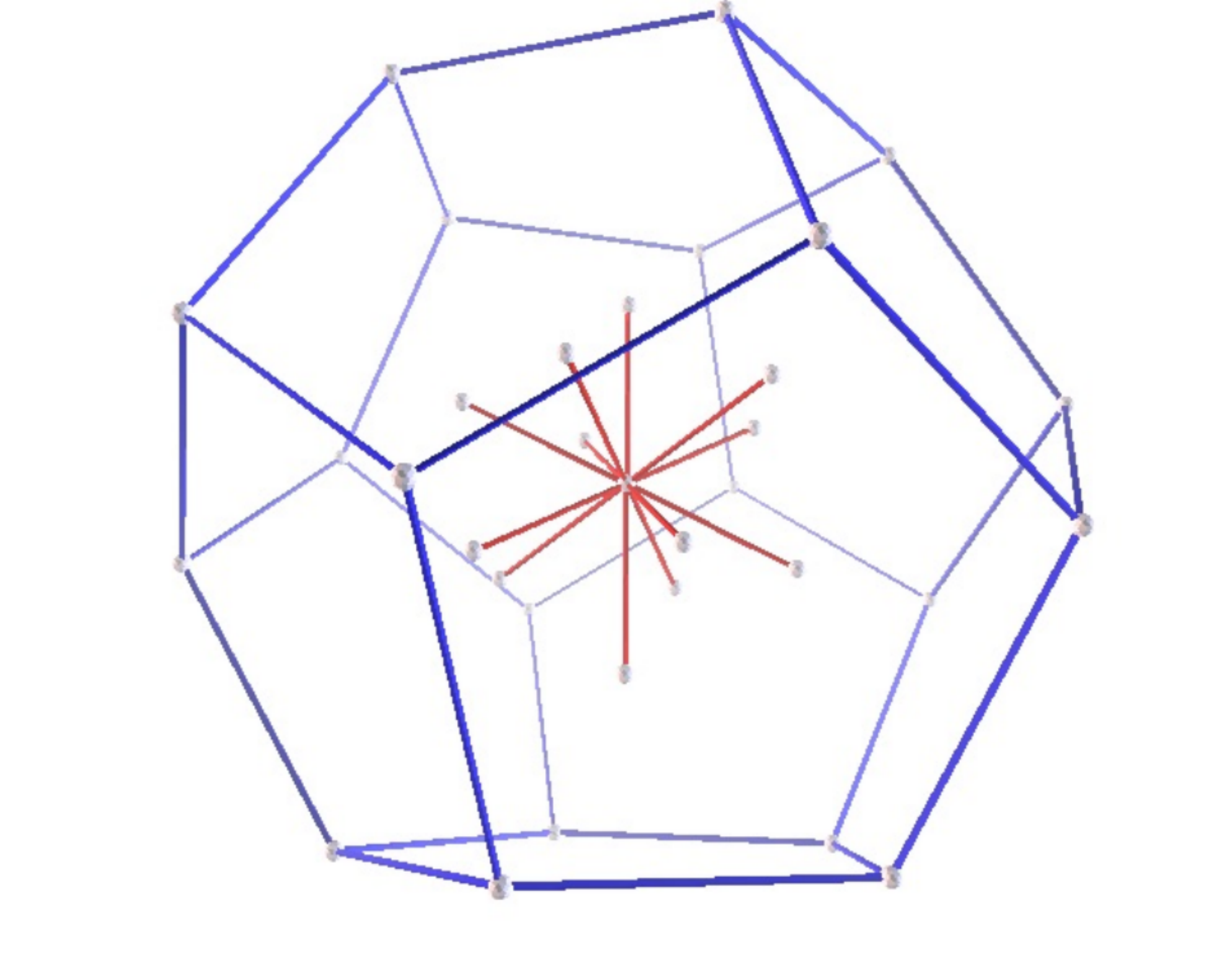}
\caption{The regular dodecahedron and its normals.}
\label{dodecahedron}
\end{center}
\end{figure}

Let us consider the regular dodecahedron $\Delta$ having vertices
$$
\begin{array}{l}
(\pm 1,\pm 1,\pm 1)\\
(0,\pm \phi,\pm \frac{1}{\phi})\\
(\pm \frac{1}{\phi} ,0,\pm \phi)\\
(\pm \phi,\pm \frac{1}{\phi},0).
\end{array}
$$
$\Sigma_\Delta$ is not rational, but it is 
quasirational with respect to the quasilattice $P\subset\mathbb{R}^3$ generated by the vectors
$$
\begin{array}{l}
Y_1=(\frac{1}{\phi},1,0)\\
Y_2=(0,\frac{1}{\phi},1)\\
Y_3=(1,0,\frac{1}{\phi})\\
Y_4=(-\frac{1}{\phi},1,0)\\
Y_5=(0,-\frac{1}{\phi},1)\\
Y_6=(1,0,-\frac{1}{\phi}),
\end{array}
$$
which is known in physics as the {\em simple icosahedral lattice}. 
In fact, the vectors $X_j=Y_j$, $X_{6+j}=-Y_j$, $j=1,\ldots,6$, are normals
for $\Delta$ (see Figure~\ref{dodecahedron}).
Notice that we have the following relations
\begin{equation}\label{relazioni}
\left\{
\begin{array}{llrrrrr}
Y_4&=&\frac{1}{\phi}Y_1&+&\frac{1}{\phi}Y_2&-&Y_3\\
Y_5&=&-Y_1&+&\frac{1}{\phi}Y_2&+&\frac{1}{\phi}Y_3\\
Y_6&=&\frac{1}{\phi}Y_1&-&Y_2&+&\frac{1}{\phi}Y_3.\\
\end{array}
\right.
\end{equation}

Let us apply the generalized toric construction to the fundamental triple
$$
(\Sigma_{\Delta}, P, \{X_1, X_2, X_3,X_4,X_5,X_6,X_7,X_8,X_9,X_{10},X_{11},X_{12}\}).
$$
The corresponding toric quasifold was first introduced in \cite{dodecahedron} and further studied,
 together with the toric spaces corresponding to all the other regular polyhedra, in
\cite{platonics}. It is given by
$$
X_{\Delta}= \frac{\bigcup_{j=1}^{20}\{\,\underline{z}\in\C^{12}\mid z_j\neq0,\;\;j\notin I_{\sigma_j}\,\}}{N_{\C}},
$$
where the $\sigma_j$'s are the maximal cones of $\Sigma_\Delta$; we identify each of them in Table~1,
 together with the corresponding vertex, fixed point, and index set $I_{\sigma_j}$.
\begin{table}
\caption{The maximal cones for the dodecahedron}
\begin{tabular}{|l|l|l|l|}
    \hline
$\sigma$ &vertex&fixed point& $I_{\sigma}$ \\ \hline
$\sigma_1$&$(-1,-1-1)$&$[0:0:0:1:1:1:1:1:1:1:1:1]$& $1,2,3$\\ \hline
$\sigma_2$&$(0,-\phi,-\frac{1}{\phi})$&$[0:0:1:0:1:1:1:1:1:1:1:1]$&$1,2,4$\\ \hline
$\sigma_3$&$(-\phi,-\frac{1}{\phi},0)$&$[0:1:0:1:1:0:1:1:1:1:1:1]$&$1,3,6$\\ \hline
$\sigma_4$&$(0,-\phi,\frac{1}{\phi})$&$[0:1:1:0:1:1:1:1:1:1:0:1]$&$1,4,11$\\ \hline
$\sigma_5$&$(-1,-1,1)$&$[0:1:1:1:1:0:1:1:1:1:0:1]$& $1,6,11$\\ \hline
$\sigma_6$&$(-\frac{1}{\phi},0,-\phi)$&$[1:0:0:1:0:1:1:1:1:1:1:1]$&$2,3,5$\\ \hline
$\sigma_7$&$(1,-1,-1)$&$[1:0:1:0:1:1:1:1:1:1:1:0]$&$2,4,12$\\ \hline
$\sigma_8$&$(\frac{1}{\phi},0,-\phi)$&$[1:0:1:1:0:1:1:1:1:1:1:0]$&$2,5,12$\\ \hline
$\sigma_9$&$(-1,1,-1)$&$[1:1:0:1:0:1:1:1:1:0:1:1]$&$3,5,10$\\ \hline
$\sigma_{10}$&$(-\phi,\frac{1}{\phi},0)$&$[1:1:0:1:1:0:1:1:1:0:1:1]$&$3,6,10$\\ \hline
$\sigma_{11}$&$(1,-1,1)$&$[1:1:1:0:1:1:1:1:0:1:0:1]$&$4,9,11$\\ \hline
$\sigma_{12}$&$(\phi,-\frac{1}{\phi},0)$&$[1:1:1:0:1:1:1:1:0:1:1:0]$&$4,9,12$\\ \hline
$\sigma_{13}$&$(0,\phi,-\frac{1}{\phi})$&$[1:1:1:1:0:1:0:1:1:0:1:1]$&$5,7,10$\\ \hline
$\sigma_{14}$&$(1,1,-1)$&$[1:1:1:1:0:1:0:1:1:1:1:0]$& $5,7,12$\\ \hline
$\sigma_{15}$&$(-1,1,1)$&$[1:1:1:1:1:0:1:0:1:0:1:1]$& $6,8,10$\\ \hline
$\sigma_{16}$&$(-\frac{1}{\phi},0,\phi$)&$[1:1:1:1:1:0:1:0:1:1:0:1]$& $6,8,11$\\ \hline
$\sigma_{17}$&$(1,1,1)$&$[1:1:1:1:1:1:0:0:0:1:1:1]$&  $7,8,9$\\ \hline
$\sigma_{18}$&$(0,\phi,\frac{1}{\phi})$&$[1:1:1:1:1:1:0:0:1:0:1:1]$& $7,8,10$\\ \hline
$\sigma_{19}$&$(\phi,\frac{1}{\phi},0)$&$[1:1:1:1:1:1:0:1:0:1:1:0]$& $7,9,12$\\ \hline
$\sigma_{20}$&$(\frac{1}{\phi},0,\phi)$&$[1:1:1:1:1:1:1:0:0:1:0:1]$& $8,9,11$\\ \hline
\end{tabular}
\label{dodecones}
\end{table}
On the other hand, $N_{\C}=\exp(\mathfrak{n}_{\C})\subset T_{\C}^{12}$, 
where $\mathfrak{n}_{\C}$ is the 
$9$--dimensional vector subspace of $\C^{12}$ having equations
$$
\left\{
\begin{array}{llrrr}
u_1&=&u_7+u_5-u_{11}&+&\frac{1}{\phi}(u_{10}-u_4+u_{12}-u_6)\\
u_2&=&u_8+u_6-u_{12}&+&\frac{1}{\phi}(u_{10}-u_4+u_{11}-u_5)\\
u_3&=&u_9+u_4-u_{10}&+&\frac{1}{\phi}(u_{11}-u_5+u_{12}-u_6).\\
\end{array}
\right.
$$
We have $20$ charts, one around each different fixed point.

Let us first consider the first two appearing in Table~\ref{dodecones}. Notice that the two index sets,
$\{1,2,3\}$ and $\{1,2,4\}$, have two indices in common. This entails that the corresponding
vertices belong to the same edge of $\Delta$.

The first chart is given by
$(\mathbb{C}^3,\Gamma_{123}, V_{123}, \eta_{123})$, where 
$$
\Gamma_{123}=\left\{\,(e^{2\pi i\phi(h+k)},e^{2\pi i\phi (h+l)},e^{2\pi i\phi (l+k) })\in (S^1)^3\mid h,k,l \in\mathbb{Z}\right\},
$$
$$V_{123}=\{[\underline{z}] \in X_{\Delta}\mid z_i\neq 0, i\neq 1,2,3\},$$
$$
\begin{array}{ccc}
\mathbb{C}^3/\Gamma_{123}&\stackrel{\eta_{123}}{\longrightarrow}& V_{123}\\
\,[z_1:z_2:z_3]&\longmapsto&\left[z_1:z_2:z_3:1:1:1:1:1:1:1:1:1\right].
\end{array}
$$
It is now useful to notice that the relations (\ref{relazioni}) can be rewritten as follows
$$
\left\{
\begin{array}{llrrrrr}
Y_3&=&\frac{1}{\phi}Y_1&+&\frac{1}{\phi}Y_2&-&Y_4\\
Y_5&=&-\frac{1}{\phi}Y_1&+&Y_2&-&\frac{1}{\phi}Y_4\\
Y_6&=&Y_1&-&\frac{1}{\phi}Y_2&-&\frac{1}{\phi}Y_4.\\
\end{array}
\right.
$$
The second chart is given by
$(\mathbb{C}^3,\Gamma_{124}, V_{124}, \eta_{124})$, where 
$$
\Gamma_{124}=\left\{\,(e^{2\pi i\phi(h+k)},e^{2\pi i\phi (h+l)},e^{2\pi i\phi (l+k) })\in (S^1)^3\mid h,k,l \in\mathbb{Z}\right\},
$$
$$V_{124}=\{[\underline{z}] \in X_{\Delta}\mid z_i\neq 0, i\neq 1,2,4\},$$
$$
\begin{array}{ccc}
\mathbb{C}^3/\Gamma_{124}&\stackrel{\eta_{124}}{\longrightarrow}& V_{124}\\
\,[w_1:w_2:w_4]&\longmapsto&\left[w_1:w_2:1:w_4:1:1:1:1:1:1:1:1\right].
\end{array}
$$
The chart transition and generalized Laurent monomials are given by
$$
\begin{array}{ccc}
\eta_{123}^{-1}(V_{123}\cap V_{124})&\stackrel{\eta_{124}^{-1}\circ\eta_{123}}
\longrightarrow& \eta_{124}^{-1}(V_{123}\cap V_{124})\\
\,[z_1:z_2:z_3]&\longmapsto&\left[z_1z_3^{\frac{1}{\phi}} : z_2z_3^{\frac{1}{\phi}}:z_3^{-1}\right].
\end{array}
$$

\begin{figure}[h]
 \begin{center}
\includegraphics[scale=0.7]{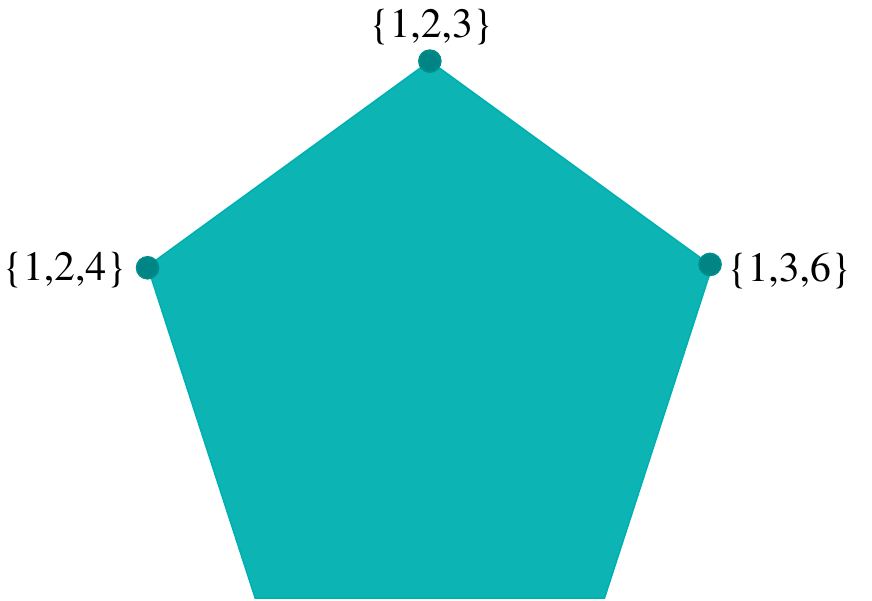}
\caption{The three vertices on their common facet.}
\label{pentagon}
\end{center}
\end{figure}

Let us now compare the second chart with the one around the third fixed point in Table~\ref{dodecones}. The 
corresponding index sets, $\{1,2,4\}$ and $\{1,3,6\}$, intersect in only one index, which entails that the 
vertices belong to the same facet, but not the same edge of the dodecahedron (see Figure~\ref{pentagon} 
for the mutual placement of all three vertices on their common facet).
Let us show how this reflects on further complexity of the Laurent monomials.
Notice first that we have
$$
\left\{
\begin{array}{llrrrrr}
Y_2&=&\frac{1}{\phi}Y_1&+&\frac{1}{\phi}Y_3&-&Y_6\\
Y_4&=&Y_1&-&\frac{1}{\phi}Y_3&- &\frac{1}{\phi}Y_6\\
Y_5&=&-\frac{1}{\phi}Y_1&+&Y_3&-&\frac{1}{\phi}Y_6.\\
\end{array}
\right.
$$
The third chart is given by $(\mathbb{C}^3,\Gamma_{136}, V_{136}, \eta_{136})$, where 
$$
\Gamma_{136}=\left\{\,(e^{2\pi i\phi(h+k)},e^{2\pi i\phi (h+l)},e^{2\pi i\phi (l+k) })\in (S^1)^3\mid h,k,l \in\mathbb{Z}\right\},
$$
$$V_{136}=\{[\underline{z}] \in X_{\Delta}\mid z_i\neq 0, i\neq 2,4,5\},$$
$$
\begin{array}{ccc}
\mathbb{C}^3/\Gamma_{136}&\stackrel{\eta_{136}}{\longrightarrow}& V_{136}\\
\,[w_1:w_3:w_6]&\longmapsto&\left[w_1:1:w_3:1:1:w_6:1:1:1:1:1:1:1\right].
\end{array}
$$
The change of charts and generalized Laurent monomials are now given by
$$
\begin{array}{ccc}
\eta_{124}^{-1}(V_{124}\cap V_{136})&\stackrel{\eta_{136}^{-1}\circ\eta_{124}}
\longrightarrow& \eta_{136}^{-1}(V_{124}\cap V_{136})\\
\,[z_1:z_2:z_4]&\longmapsto&
\left[z_1z_2^{\frac{1}{\phi}}z_4 : z_2^{\frac{1}{\phi}}z_4^{-\frac{1}{\phi}}:z_2^{-1}z_4^{-\frac{1}{\phi}}\right].
\end{array}
$$

\bigskip

{\small 

\noindent
Dipartimento di Matematica e Informatica "U. Dini", Universit\`a di Firenze \\
Viale Morgagni 67/A, 50134 Firenze, ITALY

\noindent
fiammetta.battaglia@unifi.it, elisa.prato@unifi.it}


\begin{thebibliography}{AA}

\bibitem{audin} M. Audin, The topology of torus actions on symplectic manifolds,
Progress in Mathematics \textbf{93}, Birkh\"auser, 1991.

\bibitem{cx} F. Battaglia, E. Prato, Generalized toric varieties for
simple nonrational convex polytopes, \textit{Int. Math. Res. Not.} \textbf{24} (2001), 1315--1337.

\bibitem{kite} F. Battaglia, E. Prato, The symplectic Penrose kite, \textit{Comm. Math. Phys.} \textbf{299} (2010), 577--601.

\bibitem{platonics} F. Battaglia, E. Prato, Toric geometry of the regular convex polyhedra, \textit{J. Math.} (2017), Article ID 2542796, 15 pages.

\bibitem{cut} F. Battaglia, E. Prato, Nonrational symplectic toric cuts, \textit{Internat. J. Math.} \textbf{29} (2018), 19 pages.

\bibitem{whatis} F. Battaglia, E. Prato, Nonrational polytopes and fans in toric geometry, preprint arXiv:2205.00417, (2022), 17 pages.

\bibitem{hirze} F. Battaglia, E. Prato, D. Zaffran, Hirzebruch surfaces in a one--parameter family,
\textit{Boll. Unione Mat. Ital.} {\bf 12} (2019), 293--305.

\bibitem{cox} D. Cox, J. Little, H. Schenck, Toric varieties, Graduate Studies in Mathematics
\textbf{124}, American Mathematical Society, 2011.

\bibitem{danilov} V. I. Danilov, The geometry of toric varieties, \textit{Russian Math. Surveys} \textbf{33} (1978), 97--154.

\bibitem{delzant} T. Delzant, Hamiltoniens p\'eriodiques et images convexes
de l'application moment, \textit{Bull. Soc. Math. France} \textbf{116} (1988), 315--339.

\bibitem{demazure} M. Demazure, Sous--groupes alg\'ebriques de rang maximum du groupe de Cremona,
\textit{Ann. Sci. \'Ec. Norm. Sup\'er.} \textbf{3} (1970), 507--588.

 \bibitem{di} P. Donato, P. Iglesias, Exemples de groupes différentiels~: flots irrationnels sur le tore,
\textit{C. R. Acad. Sci. Paris Sér. I Math.} \textbf{301} (1985), 127--130.

\bibitem{guillemin}V. Guillemin, Moment maps and combinatorial invariants
of Hamiltonian $T^n$--spaces, Progress in Mathematics \textbf{122}, Birkh\"auser, 1994.

\bibitem{IZP} P. Iglesias--Zemmour, E. Prato, Quasifolds, Diffeology and Noncommutative Geometry, 
\textit{J. Noncommut. Geom.} \textbf{15} (2021), 735--759.

\bibitem{KM} Y. Karshon, D. Miyamoto, Quasifold groupoids and diffeological quasifolds,
preprint arXiv:2206.14776 (2022).

\bibitem{pentaplexity} R. Penrose, Pentaplexity A Class of Non--Periodic Tilings of the Plane, \textit{Math. Intelligencer} 
\textbf{2} (1979), 32--37.

\bibitem{pcras} E. Prato, Sur une g\'en\'eralisation de la notion de V--vari\'et\'e, \textit{C. R. Acad. Sci. Paris Sér. I Math.} 
\textbf{328} (1999), 887--890.

\bibitem{p} E. Prato, Simple non--rational convex polytopes via symplectic geometry, \textit{Topology} \textbf{40} (2001), 961--975.

\bibitem{dodecahedron} E. Prato, Symplectic toric geometry and the regular dodecahedron, \textit{J. Math.} (2015), Article ID 967417, 5 pages.

\bibitem{p4} E. Prato, Toric quasifolds, \textit{Math. Intelligencer} (2022), doi:10.1007/s00283-022-10212-y.

\bibitem{ziegler} G. Ziegler, Lectures on polytopes, Graduate Texts in Mathematics \textbf{152}, Springer, 1995.

\end{thebibliography}
\end{document}